\documentclass[11pt]{article}
\usepackage{amscd}
\usepackage{amsfonts}
\usepackage{amsmath}
\usepackage{amssymb}
\usepackage{amsthm}
\usepackage{bbm}
\usepackage{CJK}
\usepackage{fancyhdr}
\usepackage{graphicx}
\usepackage{hyperref}
\usepackage{indentfirst}
\usepackage{latexsym}
\usepackage{mathrsfs}
\usepackage{xypic}
\usepackage[compress]{cite}
\newtheorem{theorem}{Theorem}[section]

\newtheorem{definition}[theorem]{Definition}
\newtheorem{proposition}[theorem]{Proposition}
\newtheorem{example}[theorem]{Example}

\newtheorem{remark}[theorem]{Remark}

\usepackage[top=1in,bottom=1in,left=1.25in,right=1.25in]{geometry}
\def\<{\langle}
\def\>{\rangle}
\def\a{\alpha}
\def\b{\beta}

\def\c{\cdot}

\def\o{\otimes}

\date{}
\begin{document}
\renewcommand{\baselinestretch}{1.2}
\renewcommand{\arraystretch}{1.0}
\title{\bf On 3-Hom-Lie-Rinehart algebras}
\author{{\bf Shuangjian Guo$^{1}$, Xiaohui Zhang$^{2}$,  Shengxiang Wang$^{3}$\footnote
        { Corresponding author(Shengxiang Wang):~~wangshengxiang@chzu.edu.cn} }\\
{\small 1. School of Mathematics and Statistics, Guizhou University of Finance and Economics} \\
{\small  Guiyang  550025, P. R. of China} \\
{\small 2.  School of Mathematical Sciences, Qufu Normal University}\\
{\small Qufu  273165, P. R. of China}\\
{\small 3.~ School of Mathematics and Finance, Chuzhou University}\\
 {\small   Chuzhou 239000,  P. R. of China}}
 \maketitle
\begin{center}
\begin{minipage}{13.cm}

{\bf \begin{center} ABSTRACT \end{center}}
We introduce the notion of 3-Hom-Lie-Rinehart algebra  and systematically describe a cohomology complex by considering
coefficient modules.  Furthermore, we consider extensions of a  3-Hom-Lie-Rinehart algebra and characterize the first cohomology space in terms of the group of automorphisms of an $A$-split abelian extension and the equivalence classes of $A$-split abelian extensions. Finally, we study formal deformations of 3-Hom-Lie-Rinehart algebras.


 \end{minipage}
 \end{center}
 \normalsize\vskip1cm

\section*{INTRODUCTION}
\def\theequation{0. \arabic{equation}}
\setcounter{equation} {0}
The notion of Lie-Rinehart algebra plays an important role in many branches of mathematics. The idea of
this notion goes back to the work of Jacobson to study certain field extensions. It is also appeared in some
different names in several areas which includes differential geometry and differential Galois theory. In
\cite{Mackenzie05}, Mackenzie provided a list of 14 different terms mentioned for this notion. Huebschmann viewed Lie-Rinehart algebras as an algebraic
counterpart of Lie algebroids defined over smoothmanifolds.  His work on several aspects of this algebra
has been developed systematically through a series of articles namely \cite{Huebschmann90,Huebschmann98, Huebschmann99, Huebschmann04} .

The notion of Hom-Lie algebras was introduced by Hartwig, Larsson, and Silvestrov in \cite{Hartwig06} as part of
a study of deformations of the Witt and the Virasoro algebras. In a Hom-Lie algebra, the Jacobi identity
is twisted by a linear map, called the Hom-Jacobi identity. Some $q$-deformations of the Witt and the
Virasoro algebras have the structure of a Hom-Lie algebra . Because of the close relation to discrete and
deformed vector fields and differential calculus \cite{H99}, many mathematicians
pay special attention to this algebraic structure.  In the sequel, many concepts and properties have been
extended to this framework of Hom-structures. The study of Hom-associative algebras, Hom-Poisson
algebras, Non-commutative Hom-Poisson algebras, Hom-Leibniz algebras, 3-Hom-Lie algerbas and most of the results anal-
ogous to the classical algebras followed in the works of Ammar et al. \cite{Ammar11}, Hartwig et al. \cite{Hartwig06}, Makhlouf
and Silvestrov \cite{Makhlouf08, Makhlouf10}, Sheng \cite{Sheng12},  Liu  et al. \cite{Liu15} and Yau \cite{Yau09}.

 Mandal and  Mishra defined modules over a Hom-Lie-Rinehart algebra and studied a cohomology with coefficients
in a left module. It presented the notion of extensions of Hom-Lie-Rinehart
algebras and deduced a characterisation of low dimensional cohomology spaces in terms of the
group of automorphisms of certain abelian extensions and the equivalence classes of those
abelian extensions in the category of Hom-Lie-Rinehart algebras in \cite{Mandal2017}. Later,  as a generalization of \cite{Castiglioni18},      they introduced a
non-abelian tensor product in the category of Hom-Lie-Rinehart algebras and  interpreted universal central extensions (and universal $\a$-central extensions) in
terms of this non-abelian tensor product in \cite{Mandal2018}.    They also explored  a relationship between Hom-Lie-Rinehart algebras
and Hom-Batalin-Vilkovisky algebras in \cite{Mandal17}.   In a sequel, they studied formal deformations of Hom-Lie-Rinehart algebras. The associated
deformation cohomology that controls deformations was constructed using multiderivations
of Hom-Lie-Rinehart algebras in \cite{Mandal18}.  Moreover, Zhang et al. stuied  crossed modules for Hom-Lie-Rinehart
algebras in \cite{zhang18}, we studied  the structures of split regular Hom-Lie Rinehart algebras in \cite{Wang19}.

Recently, Bai  et al. introduced  a class of 3-algebras which are called 3-Lie-Rinehart algebras and  discussed the basic
structures, actions and crossed modules of 3-Lie-Rinehart algebras,   also studied the derivations from 3-Lie-Rinehart algebras to 3-Lie $A$-algebras in \cite{Bai19}. Combine  \cite{Mandal2017} with \cite{Bai19}, the following questions arise naturally:
1. How do we introduce the notion of 3-Hom-Lie-Rinehart algebra?
2. How do we give the cohomology  of 3-Hom-Lie-Rinehart algebras?
3. How do we define an abelian extension of a 3-Hom-Lie-Rinehart algebra?
4.    How do we study formal deformations of 3-Hom-Lie-Rinehart algebras using multiderivations
of 3-Hom-Lie-Rinehart algebras?

The aim of this article is to answer these questions.
In Section 2,  we introduce the notion of 3-Hom-Lie-Rinehart algebras and  give some examples.
 In Section 3, we develop  the cohomology of 3-Hom-Lie-Rinehart algebras.
 In Section 4, we consider extensions of a  3-Hom-Lie-Rinehart algebra and characterize the first cohomology space in terms of the group of automorphisms of an $A$-split abelian extension and the equivalence classes of $A$-split abelian extensions.
 In Section 5,  we study formal deformations of 3-Hom-Lie-Rinehart algebras. The associated
deformation cohomology that controls deformations is constructed using multiderivations
of 3-Hom-Lie-Rinehart algebras.

\section{Preliminaries}
\def\theequation{\arabic{section}.\arabic{equation}}
\setcounter{equation} {0}

Let $R$ denote a commutative ring with unity and $\mathbb{Z}_{+}$ be the set of all non-negative integers. We will
consider all modules, algebras and their tensor products over such a ring $R$ and all linear maps to be
$R$-linear unless otherwise stated. And we recall some basic definitions and results related to our paper from \cite{Ataguema10}  and  \cite{Mandal2017}.

\begin{definition}
Given an associative commutative algebra $A$, an $A$-module $M$  and an algebra
endomorphism $\phi: A \rightarrow  A$, we call an $R$-linear map $\delta: A \rightarrow M$ a $\phi$-derivation of $A$ into  $M$ if it satisfies
the required identity:
\begin{eqnarray*}
\delta(ab)=\phi(a)\delta(b)+\phi(b)\delta(a), ~~~\mbox{for any  $a, b\in A$}.
\end{eqnarray*}
\end{definition}
Let us denote by $Der_{\phi}(A)$ the set of $\phi$-derivations of $A$ into itself.

\begin{definition}
A Hom-Lie-Rinehart algebra over $(A, \phi)$ is a tuple $(A, L,  [\c, \c], \phi, \a, \rho)$, where $A$ is an
associative commutative algebra, $L$ is an $A$-module, $[\c, \c] :  L \times L \rightarrow L$ is a skew symmetric bilinear map,
$\phi: A\rightarrow A$ is an algebra homomorphism, $\a: L\rightarrow L$ is a linear map satisfying $\a([x, y])=[\a(x), \a(y)]$,
and the $R$-map $\rho: L\rightarrow Der_{\phi}(A)$  such that following conditions hold.\\
(1) The triplet $(L,  [\c, \c], \a)$ is a Hom-Lie algebra.\\
(2) $\a(a\c x)=\phi(a)\c \a(x)$ for all $a\in A, x \in L$.\\
(3) $(\rho, \phi)$ is a representation of $(L,  [\c, \c], \a)$ on $A$.\\
(4) $\rho(a\c x)=\phi(a)\c \rho(x)$ for all $a\in A, x \in L$.\\
(5) $[x, a\c y] = \phi(a)\c [x, y] + \rho(x)(a)\a(y)$ for all $a\in A, x, y\in L$.
\end{definition}
A Hom-Lie-Rinehart algebra $(A, L,  [\c, \c], \phi, \a, \rho)$ is said to be regular if the map $\phi: A\rightarrow A$ is
an algebra automorphism and $\a: L\rightarrow L$ is a bijective map.

\begin{definition}
A 3-Hom-Lie algebra is a triple $(L, [\c, \c, \c], \a)$ consisting of a vector space $L$, a 3-ary skew-symmetric
operation $[\c, \c, \c] : L\times L\times L\rightarrow  L$ and a linear map $\a: L \rightarrow L$ satisfying the following Hom-Jacobi identity
\begin{eqnarray*}
 [\a(x), \a(y), [u, v, w]]=[[x, y, u], \alpha(v), \alpha(w)] + [\alpha(u), [x, y, v], \alpha(w)] + [\alpha(u), \alpha(v), [x, y,w]],
\end{eqnarray*}
for any $x, y,  u, v, w\in L$.
\end{definition}

\section{3-Hom-Lie-Rinehart algebras}
\def\theequation{\arabic{section}. \arabic{equation}}
\setcounter{equation} {0}

\begin{definition}
A 3-Hom-Lie-Rinehart algebra over $(A, \phi)$ is a tuple $(A, L,  [\c, \c, \c], \phi, \a, \rho)$, where $A$ is an
associative commutative algebra, $L$ is an $A$-module, $[\c, \c, \c] : L\times L \times L \rightarrow L$ is a skew symmetric trilinear map,
$\phi: A\rightarrow A$ is an algebra homomorphism, $\a: L\rightarrow L$ is a linear map satisfying $\a([x, y, z])=[\a(x), \a(y), \a(z)]$,
and the $R$-map $\rho: L\times L\rightarrow Der_{\phi}(A)$  such that following conditions hold.\\
(1) The triple $(L,  [\c, \c, \c], \a)$ is a 3-Hom-Lie algebra.\\
(2) $\a(a\c x)=\phi(a)\c \a(x)$ for all $a\in A, x \in L$.\\
(3) $(\rho, \phi)$ is a representation of $(L,  [\c, \c, \c], \a)$ on $A$.\\
(4) $\rho(a\c x, y)=\rho(x, a\c y)=\phi^2(a)\c \rho(x, y)$ for all $a\in A, x, y \in L$.\\
(5) $[x, y, a\c z] = \phi^2(a)\c [x, y, z] + \rho(x, y)(a)\a^2(z)$ for all $a\in A, x, y, z\in L$.
\end{definition}
A 3-Hom-Lie-Rinehart algebra $(A, L,  [\c, \c, \c], \phi, \a, \rho)$ is said to be regular if the map $\phi: A\rightarrow A$ is
an algebra automorphism and $\a: L\rightarrow L$ is a bijective map.

\begin{example}
A 3-Lie-Rinehart algebra $L$ over $A$ with the trilinear map $[\c, \c, \c] : L\times L \times L \rightarrow L$ and the  $R$-map $\rho: L\times L\rightarrow Der(A)$ is a 3-Hom-Lie-Rinehart algebra
$(A, L,  [\c, \c, \c], \phi, \a, \rho)$ where $\a=Id_L, \phi=Id_A$  and $\rho: L\times L\rightarrow Der_{\phi}(A)=Der(A)$.

\end{example}
\begin{example}
A 3-Hom-Lie algebra $(L,  [\c, \c, \c],  \a)$ structure over an $R$-module $L$ gives the 3-Hom-Lie-Rinehart algebra $(A, L,  [\c, \c, \c], \phi, \a, \rho)$ with $A=R$, the algebra morphism $\phi=id_{R}$ and the trivial action of $L$ on $R$.
\end{example}
\begin{example}
If we consider a 3-Lie-Rinehart algebra $L$ over $A$ along with an endomorphism
\begin{eqnarray*}
(\phi, \a):(A, L)\rightarrow (A, L)
\end{eqnarray*}
in the category of 3-Lie-Rinehart algebras,  then we get a 3-Hom-Lie-Rinehart algebra $(A, L,  [\c, \c, \c]_\a, \phi, \a, \rho_{\phi})$ as follows:\\
(1) $[x, y, z]_{\a}=\a([x, y, z])$ for any $x, y, z\in L$;\\
(2) $\rho_{\phi}(x, y)(a)=\phi(\rho(x, y)(a))$ for all $a\in A, x, y \in L$.
\end{example}

\begin{example}
Let $(L, [\c, \c, \c], \a)$ be a 3-Hom-Lie algebra over $R$ and $A$ be an associative commutative
$R$-algebra with a homomorphism $\phi: A\rightarrow A$, and $(\rho, \phi)$ be a representation of $(L, [\c, \c, \c], \a)$ on
$A$. Furthermore, the map $\rho$ is a  $R$-linear map from $L\times L\rightarrow Der_{\phi}(A)$ gives the action of $L$ on $A$ via $\phi$-derivations. Then  we have a 3-Hom-Lie-Rinehart algebra $(A, g,  [\c, \c, \c]_g, \phi, \widetilde{\a}, \widetilde{\rho}_{\phi})$, where $g=A\otimes L$. More precisely we have the following:\\
(1) The $R$-trilinear bracket $[\c, \c, \c]$ on $g$ is given by
\begin{eqnarray*}
[a\o x, b\o y, c\o z]_g&:=&\phi(abc)\o [x, y, z]+\phi(ab)\rho(x, y)(c)\a(z)+\phi(bc)\rho(y, z)(a)\a(x)\\
&&+\phi(ca)\rho(z, x)(b)\a(y),
\end{eqnarray*}
for all $x, y, z\in L$ and $a, b, c\in A$. \\
(2) The $R$-linear map $\widetilde{\a}:g\rightarrow g$ is given by
\begin{eqnarray*}
\widetilde{\a}(a\o x):=\phi(a)\o \a(x),
\end{eqnarray*}
for all $x\in L$ and $a\in A$.\\
(3) The action of $g$ on $A$ via $\phi$-derivations is given by
\begin{eqnarray*}
\widetilde{\rho}(a\o x, b\o y)(c):=\phi(ab)\rho(x, y)(c),
\end{eqnarray*}
for all $x, y\in L$ and $a, b, c\in A$.
\end{example}

\begin{example}
Let $(A, L,  [\c, \c, \c]_L, \phi, \a_L, \rho_L)$  and $(A, M,  [\c, \c, \c]_M, \phi, \a_M, \rho_M)$ be 3-Hom-Lie-Rinehart algebras over $(A, \phi)$. We consider
\begin{eqnarray*}
L\times_{Der_{\phi}A} M=\{(l, m)\in L\times M: \rho_L(l)=\rho_M(m)\}.
\end{eqnarray*}
Then $(A, L\times_{Der_{\phi}A} M, [\c, \c, \c], \phi, \a, \rho)$ is a  3-Hom-Lie-Rinehart algebra, where \\
(1)  The trilinear bracket $[\c, \c, \c]$  is given by
\begin{eqnarray*}
[(l_1,m_1),(l_2,m_2), (l_3,m_3) ]:=([l_1, l_2, l_3], [m_1, m_2, m_3]);
\end{eqnarray*}
for any $l_1,l_2,l_3\in L$ and $m_1,m_2,m_3\in M$.\\
(2) The map $\a: L\times_{Der_{\phi}A} M\rightarrow L\times_{Der_{\phi}A} M$ is given by
\begin{eqnarray*}
\a(l, m):=(\a_L(l), \a_M(m));
\end{eqnarray*}
for any $l\in L$ and $m\in M$.\\
(3) The action of $L\times_{Der_{\phi}A} M$ on $A$ is given by
\begin{eqnarray*}
\widetilde{\rho}(l_1\o m_1, l_2\o m_2)(a):=\rho_L(l_1, l_2)(a)=\rho_M(m_1, m_2)(a),
\end{eqnarray*}
for any $l_1,l_2\in L$ and $m_1,m_2\in M$ and $a\in A$.
\end{example}

\begin{example}
Let $(A, L,  [\c, \c, \c]_L, \phi, \a, \rho)$ be a 3-Hom-Lie-Rinehart algebra over $(A, \phi)$. We consider
\begin{eqnarray*}
L\times  A=\{(l, a), l\in L, a\in A\}.
\end{eqnarray*}
Then $(A, L\times A, [\c, \c, \c], \phi, \a, \rho)$ is a  3-Hom-Lie-Rinehart algebra, where \\
(1)  The trilinear bracket $[\c, \c, \c]$  is given by
\begin{eqnarray*}
&&r(x, a):=(rx, ra), (x, a)+(y, b):=(x+y, a+b), a(y, b):=(ay, ab),\\
&&[(x, a), (y, b), (z, c)]:=([x, y, z], \rho(x, y)c+\rho(y, z)a+\rho(z, x)b),
\end{eqnarray*}
for all $x, y, z\in L$,  $a, b, c\in A$ and $r\in R$.\\
(2) The map $\a: L\times A\rightarrow L\times A $ is given by
\begin{eqnarray*}
\a(x, a):=(\a(x), \phi(a));
\end{eqnarray*}
for any $x\in L$ and $a\in A$.\\
(3) The action of $L\times A$ on $A$ is given by
\begin{eqnarray*}
\widetilde{\rho}(x\o a, y\o b)(c):=\rho(x, y)(c),
\end{eqnarray*}
for any $x, y\in L$ and $a, b, c\in A$.
\end{example}
Next we define homomorphisms of 3-Hom-Lie-Rinehart algebras.
\begin{definition}
Let $(A, L,  [\c, \c, \c]_L, \phi, \a_L, \rho_L)$ and  $(B, L',  [\c, \c, \c]_{L'}, \psi, \a_{L'}, \rho_{L'})$ be  3-Hom-Lie-Rinehart algebras, then a 3-Hom-Lie-Rinehart algebra homomorphism is defined as a pair of maps $(g, f)$, where the map $g: A \rightarrow B$ is a $R$-algebra homomorphism and $f: L\rightarrow L'$ is a $R$-linear map such that
following identities hold:\\
(1) $f(a\c x) =g(a)\c f(x)$, for all $x\in L$ and $a\in A$, \\
(2) $f([x, y, z]_L)=[f(x), f(y), f(z)]_{L'}$, for all $x, y, z\in L$,\\
(3) $f(\a_L(x))=\a_{L'}(f(x))$, for all $x\in L$,\\
(4) $g(\phi(a))=\psi(g(a))$, for all $a\in A$,\\
(5) $g(\rho_L(x, y)(a))=\rho_{L'}(f(x), f(y))(g(a))$, for all $x, y\in L$ and $a\in A$.
\end{definition}
\section{Cohomology of 3-Hom-Lie-Rinehart algebras}
\def\theequation{\arabic{section}. \arabic{equation}}
\setcounter{equation} {0}
Let $A$ be an associative and commutative $R$-algebra and $\phi$ be an algebra automorphism of $A$ and $(L, \a)$
be a 3-Hom-Lie-Rinehart algebra over $(A, \phi)$.
\begin{definition} Let $M$ be an $A$-module, and $\b\in End_R(M)$. Then the pair $(M, \b)$ is a left module over a 3-Hom-Lie-Rinehart algebra  $(L, \a)$ if the following conditions hold: \\
(1) There is a linear map $\psi: L\o L\rightarrow End_R(M)$, such that the pair  $(\psi, \b)$ is a representation of $(L, [\c, \c, \c], \a)$ on $M$, \\
(2) $\b(a\c m)=\phi(a)\c \b(m)$ for all $a\in A$ and $m\in M$,\\
(3)  $\psi(a\c x, y)=\psi(x, a\c y)=\phi^2(a)\c \psi(x, y)$, for all $a\in A$ and $x, y\in L$,\\
(4) $\psi(x, y)(a\c m)=\phi^2(a)\c \psi(x, y)(m)+\rho(x, y)(a)\c \b^2(m)$, for all $x, y\in L$, $a\in A$ and $m\in M$.

\end{definition}

\begin{example}
The pair $(A, \phi)$ is a left module over $(L, \a)$.  As $(\rho, \phi)$ is a representation of $(L, [\c, \c, \c], \a)$
over $A$. Further the conditions (3) and (4) are satisfied by definition of the map $\rho$.
\end{example}

\begin{proposition}
Let $(L, \a)$ be a 3-Hom-Lie-Rinehart algebra over $(A, \phi)$ and $(M, A, \b)$ be an abelian 3-Hom-Lie $A$-algebra. Then $(M, \b)$ is a left module over $(L, \a)$ if and only if $(A, L\oplus M, \phi, \a\oplus \b, \rho_{L\oplus M})$ is a 3-Hom-Lie-Rinehart algebra with the multiplication:
\begin{eqnarray*}
&&(\a+\b)(x_1+m_1):=\a(x_1)+\b(m_1),\\
&&[x_1+m_1,x_2+m_2,x_3+m_3]:=[x_1, x_2, x_3]+\psi(x_1, x_2)m_3+\psi(x_2, x_3)m_1+\psi(x_3, x_1)m_2,\\
&&\rho_{L\oplus M}: (L\oplus M) \otimes (L\oplus M) \rightarrow Der(A), ~~\rho_{L\oplus M} (x_1+m_1, x_2+m_2):=\rho(x_1, x_2),
\end{eqnarray*}
for any $x_1, x_2, x_3\in L$ and $m_1,m_2, m_3\in M$.
\end{proposition}
{\bf Proof.}  Since $(L, \a)$ and $(M, \b)$ are $A$-modules, then $(L\oplus M, \a+\b)$ is an $A$-module and satisfies
\begin{eqnarray*}
a(x+m)=ax+am, \forall a\in A, x\in L, m\in M.
\end{eqnarray*}
 If $(M, \b)$ is a left module over $(L, \a)$. Then $(L\oplus M, \a+\b)$ is a 3-Hom-Lie algebra with above operations, and $(A, \rho_{L\oplus M})$ is a 3-Hom-Lie algebra  $L\oplus M$-module. For any $x_1, x_2, x_3\in L$,  $m_1,m_2, m_3\in M$ and $a\in A$, we have
 \begin{eqnarray*}
&&[x_1+m_1,x_2+m_2,a(x_3+m_3)]\\
&=&[x_1+m_1,x_2+m_2,ax_3]+[x_1+m_1,x_2+m_2,am_3]\\
&=&[x_1, x_2, ax_3]+\psi(x_1, x_2)am_3+\psi(x_2, ax_3)m_1+\psi(ax_3, x_1)m_2\\
&=&\phi^2(a) ([x_1, x_2, x_3]+\psi(x_1, x_2)m_3+\psi(x_2, x_3)m_1+\psi(x_3, x_1)m_2)\\
&&+\rho(x_1,x_2)(a)\alpha^{2}(x_3)+\rho(x_1,x_2)(a)\b^{2}(m_3)\\
&=& \phi^2(a)[x_1+m_1,x_2+m_2,x_3+m_3]+\rho_{L\oplus M} (x_1+m_1, x_2+m_2)(a)(\a^2+\b^{2})(x_3+m_3).
 \end{eqnarray*}
  Therefore $(A, L\oplus M, \phi, \a\oplus \b, \rho_{L\oplus M})$ is a 3-Hom-Lie-Rinehart algebra.

  Conversely, if $(A, L\oplus M, \phi, \a\oplus \b, \rho_{L\oplus M})$ is a 3-Hom-Lie-Rinehart algebra. Then for any  $x_1, x_2\in L$,  $m, m_1,m_2, m_3\in M$, we have
  \begin{eqnarray*}
&&\psi(\a(x_1), \a(x_2))(m_1, m_2, m_3)\\
&=& [\a(x_1), \a(x_2), [m_1, m_2, m_3]]\\
&=& [[x_1, x_2, m_1], \b(m_2), \b(m_3)]+[\b(m_1), [x_1, x_2, m_2], \b(m_3)]+[\b(m_1), \b(m_2), [x_1, x_2, m_3]]\\
&=& [\psi(x_1, x_2) (m_1), \b(m_2), \b(m_3)]+[\b(m_1), \psi(x_1, x_2) (m_2), \b(m_3)]+[\b(m_1), \b(m_2), \psi(x_1, x_2) (m_3)],
  \end{eqnarray*}
  and
 \begin{eqnarray*}
 \psi(x_1, x_2) (am)&=&[x_1, x_2, am]\\
 &=& \phi^{2}(a)[x_1, x_2, m]+\rho_{L\oplus M} (x_1, x_2)(a)\b^{2}(m)\\
 &=& \phi^{2}(a)\psi(x_1, x_2) (m)+\rho(x_1, x_2)(a)\b^{2}(m).
 \end{eqnarray*}
 Then $(M, \b)$ is a left module over $(L, \a)$.
  \hfill $\square$

Next we consider the $\mathbb{Z}_{+}$-graded space of $R$-modules
\begin{eqnarray*}
C^{\ast}(L; M):=\oplus_{n\geq 1}C^{n}(L; M)
\end{eqnarray*}
for 3-Hom-Lie-Rinehart algebra  $(L, \a)$  with coefficients in $(M, \b)$.  We denote by  $C^{n}(L; M)$  the space of all linear maps   $f: \wedge^2L\o ...\o \wedge^2L \wedge L\rightarrow M $ satisfying conditions below. \\
(1) $f(\a(x_1), \c \c \c, \a(x_{2n}), \a(x_{2n+1}))=\b(f(x_1, \c  \c   \c, x_{2n+1}))$,  for  all $x_i\in L,  1\leq i \leq 2n+1$,\\
(2) $f(x_1,\c \c \c, a\c x_i, \c \c \c, x_{2n+1})=\phi^{2n+1}(a)f(x_1,\c \c \c, x_i, \c \c \c, x_{2n+1})$,  for  all $x_i\in L,  1\leq i \leq 2n+1$  and $a\in A$.

Define the $R$-linear maps $\delta:  C^{n}(L; M)\rightarrow C^{n+1}(L; M)$ given by
\begin{eqnarray*}
&& \delta f(x_1,\c  \c  \c, x_{2n+1})\\
&=& (-1)^{n+1} \rho(\a^{n}(x_{2n+1}), \a^{n}(x_{2n-1}))f(x_1, \c \c \c, x_{2n-2}, x_{2n})\\
&&+ (-1)^{n+1} \rho(\a^{n}(x_{2n}), \a^{n}(x_{2n+1}))f(x_1, \c \c \c, x_{2n-1})\\
&& \sum^{n}_{k=1}(-1)^{k+1}\rho(\a^{n}(x_{2k-1}), \a^{n}(x_{2k}))f(x_1, \c \c \c, \widehat{x}_{2k-1}, \widehat{x}_{2k},\c \c \c, x_{2n+1})\\
&&+\sum^{n}_{k=1}\sum_{j=2k+1}^{2n+1}(-1)^{k}f(\a(x_1), \c \c \c, \widehat{x}_{2k-1}, \widehat{x}_{2k},\c \c \c, [x_{2k-1}, x_{2k}, x_j], \c \c  \c,  \a(x_{2n+1})).
\end{eqnarray*}
\begin{proposition}
If $f\in C^{n}(L; M)$, then $\delta f\in C^{n+1}(L; M)$ and $\delta^2=0$.
\end{proposition}
{\bf Proof.} It follows by straightforward computations. We omit details.  \hfill $\square$

By the above proposition, $(C^{\ast}(L,M), \delta)$ is a cochain complex. The resulting cohomology of the
cochain complex can be defined as the cohomology space of 3-Hom-Lie-Rinehart algebra $(L, \a)$ with
coeffcients in $(M, \b)$, and we denote this cohomology as $H^{\ast}_{HLR}(L,M)$.
\begin{definition}
Let $(L, \a)$ be a 3-Hom-Lie-Rinehart algebra  and  $(M, \b)$ be a left module over $(L, \a)$.  If $\nu \in H^0_{HLR}(L,M)$ satisfies
\begin{eqnarray*}
\rho(x, y)\nu(z)+\rho(x, z)\nu(y)+\rho(y, z)\nu(x)-\nu([x, y, z])=0,
\end{eqnarray*}
for any $x, y, z\in L$,  then $\nu$ is called a 0-cocycle associated with $\rho$.
\end{definition}
\begin{definition}
Let $(L, \a)$ be a 3-Hom-Lie-Rinehart algebra  and  $(M, \b)$ be a left module over $(L, \a)$.  If $\omega \in H^1_{HLR}(L,M)$  satisfies
\begin{eqnarray*}
&&\omega([x, u, v], \a(y), \a(z))+\omega([y, u, v], \a(z), \a(x))+\omega( \a(x), \a(y), [z, u, v])\\
&& -\omega([x, y, z], \a(u), \a(v))+\rho(\a(y), \a(z))\omega(x, u, v)+\rho(\a(z), \a(x))\omega(y, u, v)\\
&& +\rho(\a(x), \a(y))\omega(z, u, v)-\rho(\a(u), \a(v))\omega(x, y, z)=0,
\end{eqnarray*}
for any $x, y, z, u, v\in L$,  then $\omega $ is called a 1-cocycle associated with $\rho$.
\end{definition}
\section{Abelian  extensions of 3-Hom-Lie-Rinehart algebras}
\def\theequation{\arabic{section}. \arabic{equation}}
\setcounter{equation} {0}
In this section, we introduce  abelian extensions of a 3-Hom-Lie-Rinehart algebra.  We show that associated
to any abelian extension, there is a representation and a 1-cocycle. Furthermore, abelian
extensions can be classified by the first cohomology group.

\begin{definition}
The following sequence of 3-Hom-Lie-Rinehart algebras

$$\xymatrix@C=0.5cm{
  0 \ar[r] & (L'', \a'') \ar[rr]^{i} && (L', \a') \ar[rr]^{\sigma} && (L, \a)\ar[r] & 0 }$$
is a short exact sequence  if $Im(i) = Ker(\sigma), Ker(i) = 0$ and $Im(\sigma) = L$. In this case, we call $(L', \a')$ an extension of $(L, \a)$ by   a 3-Hom-Lie-Rinehart algebra $(L'', \a'')$.
 \end{definition}
An extension of 3-Hom-Lie-Rinehart algebra $(L, \a)$ is called $A$-split if we have an $A$-module map $\tau: (L, \a)\rightarrow (L', \a')$ such that\\
(1) $\sigma \circ \tau =Id_{(L, \a)}$,\\
(2) $\tau(a\c x)=a\c \tau(x)$,\\
(3) $\tau\circ \a=\a'\circ \tau$, for each $a\in A$ and $x\in L$.

Furthermore, if there exists a splitting which is also a homomorphism between 3-Hom-Lie-Rinehart algebras, we say that the extension is split.

Note that any 3-Hom-Lie-Rinehart algebra module $(M, \b)$ gives a 3-Hom-Lie-Rinehart
algebra $(A, M, [\c, \c, \c]_M, \phi, \b, \rho_M)$, with a trivial bracket and a trivial anchor map.

\begin{definition}
Let $(L, \a)$ be a 3-Hom-Lie-Rinehart algebra over $(A, \phi)$ and $(M, \b)$ be a module over $(L, \a)$. A short exact sequence
$$\xymatrix@C=0.5cm{
  0 \ar[r] & (M, \b) \ar[rr]^{i} && (L', \a') \ar[rr]^{\sigma} && (L, \a)\ar[r] & 0 }$$
is  called an abelian  extension of $(L, \a)$ by $(M, \b)$ if $(M, \b)$ is an abelian ideal of $(L, \a)$, i.e., $[\c, m, n]=0, \forall m, n\in M$.
\end{definition}

Next, we will show that the first cohomology space $H^1_{HLR}(L,M)$ of a 3-Hom-Lie-Rinehart algebra
$(L, \a)$ with coeffcients in $(M, \b)$ classifies $A$-split abelian extensions of $(L, \a)$ by $(M, \b)$.

\begin{theorem}
There is a one-to-one correspondence between the equivalent classes of $A$-split abelian
extensions of a 3-Hom-Lie-Rinehart algebra $(L, \a)$ by $(M, \b)$ and the cohomology classes in $H^1_{HLR}(L,M)$.
\end{theorem}
{\bf Proof.} Let $\omega$ be a map of the cohomology class $[\omega]\in H^1_{HLR}(L,M)$. Consider a 3-Hom-Lie-Rinehart algebra,  where the structure constraints are given as follows:\\
(1) $L'=L\oplus M$ as a direct sum of $A$-modules,\\
(2) $[x+m, y+n, z+ p]=([x, y, z]+\psi(x, y)p+\psi(y, z)m+\psi(z, x)n+\omega(x, y, z)$,\\
(3) $\a'(x+ m))=\a(x)+\b(m)$,\\
(4) $\rho'(x+m)=\rho(x)=\rho(\pi(x+ m))$,\\
for all $x, y\in L$, $m, n\in M$ and $\pi:L'\rightarrow L$ defined as $\pi(x+m)=x$.

Furthermore,
$$\xymatrix@C=0.5cm{
   & (M, \b) \ar[rr]^{i} && (L', \a') \ar[rr]^{\pi} && (L, \a)}$$
is an $A$-split abelian extension of $(L, \alpha)$ by $(M, \beta)$, where $i: M \rightarrow L'$ is defined by $i(m)= m$.

Suppose we take an another map $\omega'$ of the cohomology class $[\omega]\in H^1_{HLR}(L,M)$  and get an
extension $(L'', \a'')$ as above. Since $\omega$ and $\omega'$ represent the same cohomology class $[\omega]$, we have $\omega-\omega' = \delta\nu$
for some $\nu \in C^0(L,M)$. Then the map $F: (L', \a')\rightarrow (L'', \a'')$ defined by $F(x+m)=x+m+\nu(x)$ gives
an isomorphism of the above extensions obtained by using $\omega$ and $\omega'$ respectively. Thus for a cohomology
class in $H^1_{HLR}(L,M)$ there is a unique equivalence class of $A$-split abelian extensions of $(L, \a)$ by $(M, \b)$.

Conversely,  let
$$\xymatrix@C=0.5cm{
   & (M, \b) \ar[rr]^{i} && (L', \a') \ar[rr]^{\sigma} && (L, \a)}$$
be an $A$-split abelian extension of $(L, \alpha)$ by $(M, \beta)$. We will show that there is a 1-cocycle in $C^1(L, M)$ which is independent of a section for the map $\sigma$.

Now, we fix a section $\tau: (L, \a)\rightarrow (L', \a')$ for the map $\sigma$. Denote by
\begin{eqnarray*}
\tau(x_1, x_2)=(\tau(x_1), \tau(x_2)),
\end{eqnarray*}
and define  $\rho: L\times L \rightarrow End(M)$ by
\begin{eqnarray*}
\rho(x_1, x_2)(m):=[\tau(x_1), \tau(x_2), m],
\end{eqnarray*}
for all $x_1, x_2\in L$ and $m\in M$, it is easy to check that $\rho$ is a representation of $L$ on $M$ and does not
depend on the choice of the section $\tau$. Moreover, equivalent $A$-split abelian extensions give the same
representation of $L$ on $M$.

 Consider the map $G:  L\oplus M\rightarrow L'$ given by
\begin{eqnarray*}
G(x+m)=\tau(x)+i(m).
\end{eqnarray*}
Then $G$ is an isomorphism of $A$-modules.

Define a 1-cochain $\Omega\in C^1(L, M)$ by the following map
\begin{eqnarray*}
\Omega_{\tau} (x, y, z):= i^{-1}([\tau(x),\tau(y), \tau(z)]-\tau([x, y, z])),
\end{eqnarray*}
for all $x, y, z\in L$, we have \\
(1) $\Omega_{\tau} $ is a 3-ary skew-symmetric $R$-linear map and it satisfies  $\Omega_{\tau} (a\c x, y, z)=\phi^{2}(a)\Omega_{\tau} (x, y, z)$  for all $x, y, z\in L, a\in A$,\\
(2) $\delta(\Omega_{\tau} )=0$, which follows using 3-Hom-Jacobi identity for $(L', [\c, \c, \c]', \a')$,\\
(3) $\Omega_{\tau} \circ \a=\b\circ \Omega_{\tau} $, which follows by the relations $\tau \circ \a= \a'\circ \tau$   and $\a'\circ i = i\circ \b$.\\
Consequently, we get that $\Omega_{\tau} $ is a 1-cocycle in $C^1(L, M)$.

Finally,  we will show that for two equivalent $A$-split abelian extensions,
the associated  1-cocycles are cohomologous.

Let
$$\xymatrix@C=0.5cm{
   & (M, \b) \ar[rr]^{i'} && (L', \a') \ar[rr]^{\sigma'} && (L, \a)}$$
be another $A$-split abelian extension of $(L, \alpha)$ by $(M, \beta)$, and it is isomorphic to the extension:
$$\xymatrix@C=0.5cm{
   & (M, \b) \ar[rr]^{i} && (L', \a') \ar[rr]^{\sigma} && (L, \a)}.$$
Suppose the map $\Phi : (L', \a') \rightarrow  (L'', \a'')$ is an isomorphism of these extensions, we show that $\tau: (L, \a)\rightarrow (L', \a')$ of $\sigma$ and $\tau': (L, \a)\rightarrow (L'', \a'')$ of $\sigma'$, the respective associated cocycles $\Omega_{\tau} $ and $\Omega_{\tau'} $ are cohomologous. Consider $\tau''=\Phi\circ \tau: (L, \a)\rightarrow(L'', \a'')$ a section of $\sigma'$. Then we have $\Omega_{\tau''} =\Omega_{\tau} $. Therefore, $\Omega_{\tau} $  and  $\Omega_{\tau'} $ are cohomologous in $H^1_{HLR}(L,M)$. \hfill $\square$

In the next result we will present a characterisation of the first cohomology space $H^1_{HLR}(L,M)$ in terms of group of automorphisms of an $A$-split abelian extension.

\begin{theorem}
There is a one-to-one correspondence between the group of automorphisms of a given $A$-split
abelian extension,
$$\xymatrix@C=0.5cm{
   & (M, \b) \ar[rr]^{i} && (L', \a') \ar[rr]^{\sigma} && (L, \a)}$$
of a 3-Hom-Lie-Rinehart algebra $(L, \a)$ by $(M, \b)$ and the cohomology space $H^1_{HLR}(L,M)$.
\end{theorem}
{\bf Proof.} Similar to \cite{Mandal2017}.   \hfill $\square$

  \section{Deformations of 3-Hom-Lie-Rinehart algebras}
\def\theequation{\arabic{section}. \arabic{equation}}
\setcounter{equation} {0}
In this section, we study formal deformations of 3-Hom-Lie-Rinehart algebras. The associated
deformation cohomology that controls deformations is constructed using multiderivations
of 3-Hom-Lie-Rinehart algebras.

We recall the notion of multiderivation of degree $n$ from \cite{Mandal18}.

 \begin{definition}
 Let $M$ be an $A$-module,  $\phi: A \rightarrow  A$ be an algebra homomorphism, and
$\beta : M \rightarrow  M$ be a $\phi$-function linear map. Then a linear map
\begin{eqnarray*}
D: \wedge^{n+1}M\rightarrow M
\end{eqnarray*}
 is called a $(\phi, \b)$-multiderivation of degree $n$ (of the $A$-module $M$) if there exists a linear
map $\sigma_D: \wedge^{n}M\rightarrow Der_{\phi^n}A$ such that the following conditions are satisfied:
\begin{eqnarray*}
 &&(i)  D(\b(x_1), \b(x_2), \c \c \c, \b(x_{n+1}))=\b(D(x_1, x_2, \c \c \c, x_{n+1})),\\
 &&(ii) \sigma_D(\b(x_1), \b(x_2), \c \c \c, \b(x_{n+1}))(\phi(a))=\phi(\sigma_D(x_1, x_2, \c \c \c, x_{n+1})(a)),\\
 && (iii) \sigma_D(x_1, x_2, \c \c \c, a\c x_{n})=\phi^{n}(a)\sigma_D(x_1, x_2, \c \c \c, x_{n}),\\
 && (iv) D(x_0, x_1, \c \c \c, a\c x_{n})=\phi^{n}(a)D(x_0, x_1, \c \c \c, x_{n})+\sigma_D(x_0, x_1, \c \c \c, x_{n-1})(a)\b^{n}(x_n),
\end{eqnarray*}
for all $x_0, \c \c \c., x_n\in M$ and $a\in A$.
 The map $\sigma_D$ is called the symbol map of the $(\phi, \b)$-multiderivation $D$. Let us denote the
space of $n$-degree $(\phi, \b)$-multiderivations of $M$ by $\mathfrak{Der}_{\phi}^{n}(M, \b)$.
 \end{definition}

 By \cite{Mandal18}, we have
 \begin{eqnarray*}
\mathfrak{Der}_{\phi}^{\ast}(M, \b):=\oplus_{n\geq 0}\mathfrak{Der}_{\phi}^{n}(M, \b),
 \end{eqnarray*}
 for any $D_1\in \mathfrak{Der}_{\phi}^{p}(M, \b)$ and $D_2\in \mathfrak{Der}_{\phi}^{q}(M, \b)$, then define a bracket as follows:
 \begin{eqnarray*}
[D_1, D_2]:=(-1)^{pq}D_1\circ D_2-D_2\circ D_1,
 \end{eqnarray*}
the product $D_1\circ D_2$ is given by the expression below for any $x_0, \c \c \c, x_p, \c \c \c, x_{p+q}$,
\begin{eqnarray*}
&&(D_1 \circ D_2)(x_0, \c \c \c, x_p, \c \c \c, x_{p+q})\\
&= &\sum_{\tau\in Sh(q+1, p)}(-1)^{|\tau|}D_1(D_2(x_{\tau(0)}, . . . , x_{\tau(q)}, \b^{q}(x_{\tau(q+1)}), . . . , \b^{q}(x_{\tau(p+q)}))).
\end{eqnarray*}
Here, we denote the $(q + 1, p)$ shuffles in $S_{q+p+1}$ (the symmetric group on
the set $\{1, ¡¤ ¡¤ ¡¤ , p + q + 1\}$)  by $Sh(q + 1, p)$, and for any permutation $\tau \in S_{q+p+1}$, $|\tau|$ denotes the signature
of the permutation $(\tau)$.   Hence, the bracket $[D_1, D_2]\in \mathfrak{Der}_{\phi}^{p+q}(M, \b)$ with the symbol map $\sigma_{[D_1, D_2]}$.  Then the space of $(\phi, \b)$-multiderivations of $M$ has a graded Lie algebra structure.

Next we describe a 3-Hom-Lie-Rinehart algebra structures in terms of the
graded Lie algebra obtained above.
 \begin{proposition}
 Let $L$ be an $A$-module and $\a: L \rightarrow L$ be a $\phi$-function linear map. Then
there is a one-to-one correspondence between 3-Hom-Lie-Rinehart algebra structures on the
pair $(L, \a)$ and elements $\mathrm{m}\in \mathfrak{Der}_{\phi}^{2}(L, \a)$ satisfying $\mathrm{m}\circ\mathrm{m} = 0$.
 \end{proposition}
 {\bf Proof.} Let $(L, \a)$ be a 3-Hom-Lie-Rinehart algebra over $(A, \phi)$. Define a trilinear map $\mathrm{m}: L\times L \times L\rightarrow L$ by $\mathrm{m}(x, y, z):=[x, y, z]$, for any $x, y, z\in L$. By definition for any $x, y, z\in L$ and $a\in A$, we have
 \begin{eqnarray}
 \mathrm{m}(x, y, a\c z)=\phi^{2}(a) \mathrm{m}(x, y, z)+\rho(x, y)(a)\a^{2}(z).
 \end{eqnarray}
 It follows that $ \mathrm{m}$ is a  2-degree $(\phi, \a)$-derivation of the $A$-module $L$, i.e. $\mathrm{m}\in \mathfrak{Der}_{\phi}^{2}(L, \a)$ with symbol $\sigma_{\mathrm{m}}=\rho: L\times L\rightarrow Der_{\phi}(A)$. Furthermore, from the definition of the graded Lie bracket, we calculate
 \begin{eqnarray*}
&&\mathrm{m}\circ \mathrm{m}(x, y, u, v, w)\\
&=&[\a(x), \a(y), [u, v, w]]-[\alpha(v), \alpha(w), [x, y, u]] +[\alpha(u),  \alpha(w),  [x, y, v]] - [\alpha(u), \alpha(v), [x, y,w]]\\
&=& 0.
 \end{eqnarray*}

 Conversely, suppose  that $\mathrm{m}\in \mathfrak{Der}_{\phi}^{2}(L, \a)$ satisfies $\mathrm{m}\circ\mathrm{m} = 0$. Let us define a bracket $[\c, \c, \c]: L\times L \times L\rightarrow L$ as follows:
 \begin{eqnarray*}
[x, y, z]:=\mathrm{m}(x, y, z), ~~\mbox{for any~} x, y, z\in L.
 \end{eqnarray*}
  Also define a linear map $\rho:=\sigma_{\mathrm{m}}: L\times L\rightarrow Der_{\phi}(A)$. It is easy  to verify that $(A, L,  [\c, \c, \c], \phi, \a, \rho)$ is a  3-Hom-Lie-Rinehart algebra.  \hfill $\square$

Let us define a cochain complex $(C^{\ast}_{def}(L, \a), \delta)$, where
\begin{eqnarray*}
C^{\ast}_{def}(L, \a):=\oplus_{n\geq 1}C^{n}_{def}(L, \a), ~~\mbox{and}~~ C^{\ast}_{def}(L, \a):=\mathfrak{Der}_{\phi}^{n-1}(L, \a).
\end{eqnarray*}
Define the differential
$
\delta: C^{n}_{def}(L, \a)\rightarrow C^{n+1}_{def}(L, \a)
$
by
\begin{eqnarray*}
\delta(D)=[\mathrm{m},  D],
\end{eqnarray*}
for any $D\in C^{n}_{def}(L, \a)$.  In particular, for any $x_1, \c \c  \c, x_{2n+1}\in L$ and $D\in \mathfrak{Der}_{\phi}^{n-1}(L, \a)$, the
coboundary expression is given as follows.
\begin{eqnarray*}
&& \delta (D)(x_1,\c  \c  \c, x_{2n+1})\\
&=& (-1)^{n+1} \mathrm{m}(\a^{n}(x_{2n+1}), \a^{n}(x_{2n-1}), D(x_1, \c \c \c, x_{2n-2}, x_{2n}))\\
&&+ (-1)^{n+1} \mathrm{m}(\a^{n}(x_{2n}), \a^{n}(x_{2n+1}), D(x_1, \c \c \c, x_{2n-1}))\\
&& \sum^{n}_{k=1}(-1)^{k+1}\mathrm{m}(\a^{n}(x_{2k-1}), \a^{n}(x_{2k}), D(x_1, \c \c \c, \widehat{x}_{2k-1}, \widehat{x}_{2k},\c \c \c, x_{2n+1}))\\
&&+\sum^{n}_{k=1}\sum_{j=2k+1}^{2n+1}(-1)^{k}D(\a(x_1), \c \c \c, \widehat{x}_{2k-1}, \widehat{x}_{2k},\c \c \c, \mathrm{m}(x_{2k-1}, x_{2k}, x_j), \c \c  \c,  \a(x_{2n+1})).
\end{eqnarray*}
Note $\mathrm{m}\in \mathfrak{Der}_{\phi}^{2}(L, \a)$  satisfies $\mathrm{m}\circ\mathrm{m} = 0$, therefore  $\delta^2=0$.

We now consider the 3-Hom-Lie-Rinehart algebra structure on
$(L, \alpha)$ over $(A,\phi)$ as an element $\mathrm{m}\in \mathfrak{Der}_{\phi}^{2}(L, \a)$ satisfying $\mathrm{m}\circ\mathrm{m} = 0$. Here, we denote
by $R[[t]]$ the space of formal power series ring with parameter $t$.

\begin{definition}
A deformation of a 3-Hom-Lie-Rinehart algebra $(L, \a)$ over $(A, \phi)$, given
by $\mathrm{m}\in \mathfrak{Der}_{\phi}^{2}(L, \a)$, is a $R[[t]]$-bilinear map
\begin{eqnarray*}
\mathrm{m}_t: L[[t]]\times L[[t]]\times L[[t]]\rightarrow L[[t]], ~~~\mathrm{m}_t(x, y, z)=\sum_{i\geq 0}t^im_i(x, y, z),
\end{eqnarray*}
with $m_0=\mathrm{m}$ and $m_i\in \mathfrak{Der}_{\phi}^{2}(L, \a)$ for $i\geq 0$,    satisfying $\mathrm{m}_t\circ \mathrm{m}_t=0$,  and maps $\phi_t : A[[t]]\rightarrow A[[t]]$ and $\alpha_{L_t}: L[[t]] \rightarrow L[[t]]$ are extensions of the maps $\phi$ and $\a$ with $t$.
\end{definition}

Let $\mathrm{m}_t$ be a deformation of $\mathrm{m}$. Then
 \begin{eqnarray*}
&&\mathrm{m}_t(\a(x), \a(y),\mathrm{m}_t (u, v, w))-\mathrm{m}_t(\alpha(v), \alpha(w), \mathrm{m}_t(x, y, u)) +\mathrm{m}_t(\alpha(u),  \alpha(w),  \mathrm{m}_t(x, y, v)) \\
&&- \mathrm{m}_t(\alpha(u), \alpha(v), \mathrm{m}_t(x, y,w))=0.
 \end{eqnarray*}
Comparing the coefficients of $t^n$, $n\geq 0$, we get the following equation:
 \begin{eqnarray*}
&&\sum_{i, j=0}^nm_i(\a(x), \a(y),m_j (u, v, w))-m_i(\alpha(v), \alpha(w), m_j(x, y, u)) +m_i(\alpha(u),  \alpha(w),  m_j(x, y, v)) \\
&&- m_i(\alpha(u), \alpha(v), m_j(x, y,w))=0.
 \end{eqnarray*}

 \begin{remark}
 For $n=1$, then we have $[\mathrm{m}, m_1]=\delta(m_1)=0$, i.e. $m_1$ is a 3-cocycle.

 \end{remark}

  \begin{definition}
  The 3-cochain $m_1$ is called the infinitesimal of the deformation $\mathrm{m}_t$. More
generally, if $m_i = 0$ for $1 \leq i \leq (n-1)$ and $m_n$ is non zero cochain, then $m_n$ is called the
$n$-infinitesimal of the deformation $\mathrm{m}_t$.

 \end{definition}
 \begin{definition}
 Two deformations $\mathrm{m}_t$ and $\mathrm{\widetilde{m}}_t$ are said to be equivalent if there exists a formal
automorphism
 \begin{eqnarray*}
\Phi_t : L[[t]]\rightarrow L[[t]],~~~\Phi_t=id_L+\sum_{i\geq 1}t^i\phi_i,
 \end{eqnarray*}
 where $\phi_i: L\rightarrow L$ is a $R$-linear map such that
 \begin{eqnarray*}
 \phi_i\circ \a=\a\circ \phi_i,~~~  \mathrm{\widetilde{m}}(x, y, z)=\Phi_t^{-1}\mathrm{m}_t(\Phi_t(x), \Phi_t(y), \Phi_t(z)).
 \end{eqnarray*}
 \end{definition}

  \begin{definition}
A deformation is called trivial if it is equivalent to the deformation $m_0 = \mathrm{m}$.
 \end{definition}
 \begin{theorem}
 The cohomology class of the infinitesimal of a deformation $\mathrm{m}_t$ is determined
by the equivalence class of $\mathrm{m}_t$.
 \end{theorem}
 {\bf Proof.} Straightforward.   \hfill $\square$
  \begin{definition}
  A 3-Hom-Lie-Rinehart algebra is said to be rigid if  every deformation
of it is trivial.
   \end{definition}

  \begin{theorem}

  A non-trivial deformation of a 3-Hom-Lie-Rinehart algebra is equivalent to
a deformation whose $n$-infinitesimal cochain is not a coboundary for some $n\geq 1$.

  \end{theorem}
{\bf Proof.} Let $\mathrm{m}_t$ be a deformation of 3-Hom-Lie-Rinehart algebra with n-infinitesimal $m_n$ for
some $n\geq 1$. Assume that there exists a 3-cochain $\phi\in C^{2}_{def}(L, \a)$ with $\delta (\phi=m_n)$. Set
\begin{eqnarray*}
\Phi_t=id_L+\phi t^n~~~\mbox{~and~}  \mathrm{\widetilde{m}}=\Phi_t^{-1}\circ \mathrm{m}_t\circ\Phi_t
\end{eqnarray*}
Comparing the coefficients of $t^n$, we get the following equation:
\begin{eqnarray*}
\widetilde{m}_n-m_n=-[\widetilde{m}, m_n]=-\delta(\phi).
\end{eqnarray*}
So $\widetilde{m}_n=0$.  Then the deformation whose $n$-infinitesimal is not a coboundary for some $n\geq 1$.   \hfill $\square$
 \begin{center}
 {\bf ACKNOWLEDGEMENT}
 \end{center}

 The paper is supported by  the NSF of China (No. 11761017),  the Youth Project for Natural Science Foundation of Guizhou provincial department of education (No. KY[2018]155) and
 the Anhui Provincial Natural Science Foundation (Nos. 1908085MA03 and 1808085MA14).

\end{document}